\documentclass[11pt,a4paper]{article}
\usepackage[T1]{fontenc}
\usepackage{lmodern}
\usepackage[margin=29mm,headheight=14pt]{geometry}
\usepackage{amsmath,amssymb,amsthm,mathtools}
\usepackage{mathrsfs}
\usepackage{aliascnt}
\usepackage{microtype}
\usepackage{needspace}
\usepackage{booktabs,longtable,array}
\usepackage[dvipsnames]{xcolor}
\definecolor{InternalLink}{HTML}{0057B8}
\definecolor{CitationLink}{HTML}{1B7837}
\definecolor{URLLink}{HTML}{A33A00}
\usepackage{enumitem}
\setlist[enumerate]{label=\textup{(\roman*)},leftmargin=2em,itemsep=3pt,topsep=5pt}
\usepackage{hyperref}
\hypersetup{colorlinks=true,linkcolor=InternalLink,citecolor=CitationLink,urlcolor=URLLink,linktoc=all,pdfborder={0 0 0},pdftitle={One complex metric and geometric adelic nefness},pdfsubject={Adelic line bundles, complex semipositivity, and geometric nefness}}
\usepackage[nameinlink,noabbrev]{cleveref}
\usepackage{bookmark}
\usepackage{fancyhdr}
\newtheorem{theorem}{Theorem}[section]
\theoremstyle{plain}
\newaliascnt{proposition}{theorem}
\newtheorem{proposition}[proposition]{Proposition}
\aliascntresetthe{proposition}
\crefname{proposition}{proposition}{propositions}
\Crefname{proposition}{Proposition}{Propositions}
\newaliascnt{lemma}{theorem}
\newtheorem{lemma}[lemma]{Lemma}
\aliascntresetthe{lemma}
\crefname{lemma}{lemma}{lemmas}
\Crefname{lemma}{Lemma}{Lemmas}
\newaliascnt{corollary}{theorem}

\aliascntresetthe{corollary}
\crefname{corollary}{corollary}{corollaries}
\Crefname{corollary}{Corollary}{Corollaries}
\theoremstyle{definition}
\newaliascnt{definition}{theorem}
\newtheorem{definition}[definition]{Definition}
\aliascntresetthe{definition}
\crefname{definition}{definition}{definitions}
\Crefname{definition}{Definition}{Definitions}
\theoremstyle{remark}
\newaliascnt{remark}{theorem}
\newtheorem{remark}[remark]{Remark}
\aliascntresetthe{remark}
\crefname{remark}{remark}{remarks}
\Crefname{remark}{Remark}{Remarks}
\numberwithin{equation}{section}
\newcommand{\OO}{\mathcal O}
\newcommand{\LL}{\mathcal L}
\newcommand{\AdPic}{\widehat{\Pic}}
\newcommand{\II}{\mathcal I}
\newcommand{\JJ}{\mathcal J}
\newcommand{\QQ}{\mathbb Q}
\newcommand{\CC}{\mathbb C}
\newcommand{\ZZ}{\mathbb Z}

\newcommand{\an}{\mathrm{an}}
\newcommand{\geom}{\mathrm{geom}}
\newcommand{\GS}{\mathrm{GS}}

\DeclareMathOperator{\Pic}{Pic}

\DeclareMathOperator{\divisor}{div}

\DeclareMathOperator{\Supp}{Supp}
\DeclareMathOperator{\im}{im}

\newcommand{\bnorm}[1]{\left\lVert#1\right\rVert_{D_0}}
\begin{document}

\title{One complex metric and geometric adelic nefness}
\author{Kaicheng Dong}
\date{}
\maketitle

\begin{abstract}
We prove that an arithmetic adelic line bundle on a smooth
quasi-projective variety over a number field has nef geometric
part if its prescribed continuous metric at one complex embedding
has semipositive first Chern current. The proof constructs
semiample approximations to every positive rational ample
perturbation by boundary extension, multiplier ideals, and Nadel
vanishing, followed by specialization and Galois descent.
This gives geometric nefness of Gao--Zhang's canonical
Gross--Schoen adelic line bundle.
\end{abstract}


\begingroup
\small
\tableofcontents
\endgroup
\clearpage

\section{Introduction}\label{sec:introduction}
A line bundle on a smooth projective complex variety is nef
if it admits a continuous Hermitian metric with semipositive
curvature. For a quasi-projective variety, the corresponding
question depends on the extension of the line bundle at the boundary.

In this paper, a \emph{compactified $\QQ$-line bundle} on a smooth
quasi-projective variety $U$ is an equivalence class of
boundary-Cauchy sequences of $\QQ$-line bundles $P_i$ on normal
projective models $U\subset X_i$, with specified identifications
$\alpha_i:L_U\xrightarrow{\sim}P_i|_U$ for a fixed $\QQ$-line bundle
$L_U$. On a common model, the differences $P_i-P_j$ are boundary
divisors, measured by upper and lower bounds in multiples of a
fixed effective boundary divisor. Two sequences are equivalent
if their difference tends to zero in this boundary norm, after
identifying their underlying line bundles on $U$.
These are Yuan--Zhang's geometric adelic line bundles \cite{YZ};
see \Cref{def:completion} for details.

We prove that a continuous semipositive metric at one complex
place implies geometric nefness when the metric and the
compactified line bundle are limits of the same model data.
This compatibility holds for arithmetic adelic line bundles.
Yuan \cite[Section~3.3.5]{Yuan} and Gao--Zhang
\cite[Section~4.1]{GZ} each give a definition of a
Gross--Schoen adelic line bundle. Applied to Gao--Zhang's
line bundle, our result answers the geometric nefness
question in \cite[Remark~5.2]{GZ}.

\subsection{Complex semipositivity and geometric nefness}
We follow Yuan--Zhang's definitions of projective models and the
boundary topology \cite[Sections~2.4--2.5.1]{YZ}, and of strongly
nef and nef adelic line bundles \cite[Definition~2.5.2]{YZ};
see also \Cref{sec:preliminaries}.

Let $U$ be a smooth complex quasi-projective variety. Fix a smooth projective
compactification $U\subset X_0$ and an effective integral Cartier
divisor $D_0$ with $\operatorname{Supp}(D_0)=X_0\setminus U$.
We write $\|\cdot\|_{D_0}$ for the corresponding boundary norm.
Comparisons of model line bundles are understood with respect to
their specified identifications over $U$.

For a continuous Hermitian metric $h$ on a line bundle $L_U$ and a
local frame $e$, write $\varphi_{h,e}=-\log\|e\|_h$.
We say that $h$ has \emph{semipositive curvature} if
$c_1(L_U,h)=dd^c\varphi_{h,e}\geq0$ in the sense of currents;
our normalization is recalled in \Cref{def:weights}.
When $D_0\neq0$, choose a smooth metric on $\OO(D_0)$, normalized
so that its boundary Green function
$g_0=-\log\|s_{D_0}\|\geq1$ on $U$.
If $U$ is projective, take $X_0=U$, $D_0=0$, and
$\LL=L_U$; in this case, the boundary conditions below
are omitted.

\paragraph{Compactified metrics.}
Let $\LL$ be a compactified $\QQ$-line bundle on $U$ with underlying
line bundle $L_U$.
A \emph{compactified metric} on $\LL$ is a continuous Hermitian
metric $h$ on $L_U$ with the following property.
For every rational $t>0$, there exist a projective modification
$p_t:X_t\to X_0$, with $X_t$ normal and $p_t$ an isomorphism over $U$,
a $\QQ$-line bundle $P_t$ on $X_t$, a specified isomorphism
$\alpha_t:L_U\xrightarrow{\sim}P_t|_U$, and a continuous Hermitian
metric $h_t$ on $P_t$, such that
\begin{equation}\label{eq:simultaneous-approximation}
 \bnorm{P_t-\LL}\leq t,
 \qquad
 \bigl|\varphi_{\alpha_t^*h_t,e}-\varphi_{h,e}\bigr|
     \leq t g_0\quad\text{on }U.
\end{equation}
Both comparisons use the same specified isomorphism $\alpha_t$.
For projective $U$, every continuous Hermitian metric on $L_U$
is a compactified metric. See \Cref{def:compactified-metric}.

\Needspace{7\baselineskip}
\begin{theorem}[Complex semipositivity and geometric nefness]
\label{thm:main-nefness}
Let $U$ be a smooth irreducible complex quasi-projective variety
of positive dimension, and let $\LL$ be a compactified
$\QQ$-line bundle on $U$ with underlying line bundle $L_U$.
If $\LL$ admits a compactified metric $h$ satisfying
\[
 c_1(L_U,h)\geq0,
\]
then $\LL$ is nef.
\end{theorem}

The prescribed complex metric of a Yuan--Zhang arithmetic adelic
line bundle is a compactified metric on its geometric part.
This follows from the arithmetic Cauchy condition
\cite[Section~2.5.1]{YZ}; see \Cref{prop:adelic-approximation}.
A continuous semipositive metric chosen independently of $\LL$
need not satisfy this compatibility or imply nefness of the
prescribed compactified extension.

We sketch the proof in \Cref{sec:proof-sketch} and give the full
argument in \Cref{sec:main-proof}.

\subsection{The theorem over a number field}
Write $\AdPic(U/K)_{\QQ}$ and $\AdPic(U/\ZZ)_{\QQ}$ for the groups of compactified and arithmetic adelic $\QQ$-line bundles, respectively. For an arithmetic adelic line bundle $\overline{\LL}$, denote its geometric part by $\LL^{\geom}$. It is obtained by taking the generic fibers of a model presentation, with their open isomorphisms; it is a compactified line bundle over $K$. 

\Needspace{7\baselineskip}
\begin{theorem}[Geometric nefness over a number field]\label{thm:number-field}
Let $K$ be a number field, let $U/K$ be a smooth geometrically
integral quasi-projective variety of positive dimension, and
let $\LL\in\AdPic(U/K)_{\QQ}$ have underlying line bundle $L_U$.
If, for an embedding $\sigma:K\hookrightarrow\CC$, $\LL_\sigma$
admits a compactified metric $h$ satisfying
\[
 c_1(L_{U,\sigma},h)\geq0,
\]
then $\LL$ is nef.

In particular, the geometric part of an arithmetic adelic line
bundle is nef if its prescribed continuous metric at one complex
embedding has semipositive first Chern current.
\end{theorem}

The proof in \Cref{sec:arithmetic} specializes generating sections
and boundary comparison maps to a finite extension of $K$, then
descends their tensor products over the Galois conjugates.
Division by the extension degree preserves the boundary error.

\subsection{Gross--Schoen line bundles}
For a smooth projective curve $C$ of genus $g\geq2$ and a degree-one class $\xi$ with $(2g-2)\xi=\omega_C$, let $\Delta_\xi$ denote the canonical Gross--Schoen modified diagonal cycle on $C^3$.

For a smooth family $\pi:\mathcal C\to S$, Yuan constructs the admissible adelic dualizing bundle $\overline\omega_a$ and the adelic divisor $\Phi_S$ associated with the local $\varphi$-invariant. His Gross--Schoen line bundle is
\[
 \overline G^{\,Y}_{\mathcal C/S}
   =\frac{2g+1}{2g-2}
       \langle\overline\omega_a,\overline\omega_a\rangle_\pi
      -\overline{\OO}(\Phi_S).
\]
Its point heights are the canonical Gross--Schoen heights. Yuan asks whether this adelic line bundle is nef and, when the family has maximal variation, whether it is big \cite[Section~3.3.5]{Yuan}. Here nefness and bigness are arithmetic; see \Cref{def:yuan-GS} for the normalization.

Gao--Zhang construct a Gross--Schoen adelic line bundle on the parameter space of canonically centered curves by a Deligne pairing. Let $M$ be a smooth geometrically integral fine-level parameter variety over a number field $K$ for pairs $(C,\xi)$ with $g(C)=g\geq3$ and $(2g-2)\xi=\omega_C$. To obtain $M$, start with a fine moduli scheme of smooth genus-$g$ curves with level structure of order at least $3$, take the finite \'etale cover parametrizing classes $\xi$ in the relative Picard scheme satisfying $(2g-2)\xi=\omega$, and choose a geometrically integral component over a number field of definition. Let $\pi:\mathcal C\to M$ be the universal curve, $J=\operatorname{Jac}(\mathcal C/M)$, and $i_\xi:\mathcal C\to J$ the Abel--Jacobi map. If $\overline{\mathcal P}$ is their canonically metrized Poincar\'e line bundle, set
\begin{equation}\label{eq:GS-introduction}
 \overline{\mathcal Q}=(i_\xi,i_\xi)^*\overline{\mathcal P},
 \qquad
 \overline G_g
   =\langle\overline{\mathcal Q},
            \overline{\mathcal Q},
            \overline{\mathcal Q}\rangle_{\mathcal C\times_M\mathcal C/M}.
\end{equation}
We use the normalization of \cite[Section~4.1, equation~(4.2)]{GZ}, specified in \Cref{def:GZ-GS}. Gao--Zhang prove that $\overline G_g$ is integrable, that its point height is $\langle\Delta_\xi,\Delta_\xi\rangle_{\mathrm{BB}}$, and that its prescribed complex metric has curvature
\[
 c_1(\overline G_{g,\sigma})=\beta_{\GS}\geq0
\]
\cite[Theorems~4.1--4.2]{GZ}. Their Remark~5.2 asks for nefness of the geometric part, denoted $\widetilde{\mathscr L}$ there. Applying \Cref{thm:number-field} gives the following result.

\Needspace{7\baselineskip}
\begin{theorem}[Geometric nefness of the Gross--Schoen line bundle]\label{thm:GS}
Let $g\geq3$, let $K$ be a number field, and let $M/K$ be a
smooth geometrically integral fine-level parameter variety
for canonically centered curves as above.
The geometric part $G_g^{\geom}$ of the arithmetic adelic
$\QQ$-line bundle $\overline G_g$ defined by
\eqref{eq:GS-introduction} is nef.
\end{theorem}

\subsection{Sketch of the proof of Theorem~\ref{thm:main-nefness}}\label{sec:proof-sketch}

We adapt the multiplier-ideal argument of
Botero--Burgos Gil--Holmes--de Jong \cite[Section~5]{BBGH}.
The boundary estimates below give convergence to the prescribed
compactified line bundle while retaining its restriction over $U$.

\paragraph{Extension across the boundary.}
Fix an ample $\QQ$-line bundle $A$ on $X_0$ and a rational
$\varepsilon>0$, and choose a sufficiently small rational $t>0$.
Choose approximation data
$(p_t:X_t\to X_0,P_t,\alpha_t,h_t)$
as in \eqref{eq:simultaneous-approximation}, with $X_t$ smooth.
Set $D_t=p_t^*D_0$ and
$B_t=P_t+t D_t$.
Using the specified identification over $U$ and the
canonical section of the boundary divisor, we transport
the metric $h$ to $B_t|_U$.
Its local weights $w_t$ satisfy
\[
2t\log|f_{D_t}|-C
\leq w_t\leq C,
\]
where $f_{D_t}$ is a local equation of $D_t$
and $C$ is a constant on each sufficiently small chart.
The upper bound allows the metric to extend across
the boundary as a singular semipositive metric.
Moreover, the boundary approximation gives
$0\leq B_t-\LL\leq2t D_0$ in the
completed boundary order; see \Cref{prop:small-extension}.

\paragraph{Multiplier ideals.}
The lower bound for $w_t$ gives, for every sufficiently
divisible positive integer $m$, the inclusions
\[
\OO(-2mt D_t)
\subseteq\JJ(mw_t)\subseteq\OO.
\]
By coherence and GAGA, the multiplier ideal algebraizes
to a coherent ideal $\II_{t,m}$ on $X_t$.
Let $\rho:Y\to X_t$ be a principalization of this
ideal, with image ideal
$\II_{t,m}\OO_Y=\OO_Y(-F)$.
The preceding inclusions yield
\[
0\leq\frac{F}{m}\leq2t\rho^*D_t.
\]
This is an inequality of $\QQ$-Cartier divisors, including exceptional
multiplicities; see \Cref{prop:ideal-divisor}.

\paragraph{Global generation.}
Choose an ample $\QQ$-line bundle $A_t$ on $X_t$
with a specified identification $A_t|_U\simeq A|_U$
and satisfying
\[
-\lambda D_t
\leq A_t-p_t^*A
\leq\lambda D_t
\]
for an arbitrarily small rational $\lambda>0$.
Equip $A_t$ with a smooth Hermitian metric of
strictly positive curvature.
The induced singular metric on
$B_t+\varepsilon A_t$ then has strictly
positive curvature.

Nadel vanishing and Castelnuovo--Mumford regularity
\cite[Theorems~5.9--5.10 and Corollary~5.11]{BBGH} imply that
\[
\OO\bigl(m(B_t+\varepsilon A_t)\bigr)
\otimes\II_{t,m}
\]
is globally generated for all sufficiently large
divisible integers $m$.
Here $X_t$ and $A_t$ are fixed before $m$ is chosen;
see \Cref{lem:exact-generation}.

\paragraph{Semiample approximation.}
Choose a sufficiently large divisible $m$ and principalize
$\II_{t,m}$ as above.
Pulling back the generating sections and using the
principalized image ideal gives a semiample
$\QQ$-line bundle
\[
N_{t,\varepsilon,m}
=\rho^*(B_t+\varepsilon A_t)-\frac{F}{m}
\]
on $Y$, with the specified restriction
$N_{t,\varepsilon,m}|_U\simeq
L_U+\varepsilon A|_U$.

For $\lambda\leq t/(1+\varepsilon)$, the boundary estimates give
\[
\bigl\|N_{t,\varepsilon,m}
-(\LL+\varepsilon A)\bigr\|_{D_0}
\leq2t+\varepsilon\lambda
\leq3t.
\]
Letting $t\to0$ shows that $\LL+\varepsilon A$ is strongly nef.
Taking $\varepsilon=1/a$ and multiplying by $a$ gives strong
nefness of $a\LL+A$ for every integer $a\geq1$.
Since $A$ is fixed and strongly nef, $\LL$ is nef.

\begin{remark}[Semiample approximation]\label{rem:semiample-approximation}
The proof also gives semiample approximations with arbitrary
boundary error. Under the hypotheses of \Cref{thm:main-nefness},
fix an ample $\QQ$-line bundle $A$ on $X_0$.
For every rational $\varepsilon>0$ and every real $\delta>0$,
there exist a projective modification $q:Y\to X_0$, with $Y$
smooth and $q$ an isomorphism over $U$, and a semiample
$\QQ$-line bundle $N$ on $Y$, together with a specified isomorphism
\[
 L_U+\varepsilon A|_U\xrightarrow{\sim}N|_U,
\]
such that
\begin{equation}\label{eq:main-approximation}
 \bnorm{N-(\LL+\varepsilon A)}<\delta.
\end{equation}
Consequently, $\LL+\varepsilon A$ is strongly nef for every
rational $\varepsilon>0$, and $a\LL+A$ is strongly nef for
every integer $a\geq1$, with the same fixed ample line bundle $A$.
If $U$ is projective, one may take $Y=X_0=U$ and
$N=L_U+\varepsilon A$, with no boundary estimate needed.
\end{remark}

\section*{Acknowledgements}
\addcontentsline{toc}{section}{Acknowledgements}
I am grateful to Zhelun Chen for introducing me to this problem during a
conference at Bielefeld University in March 2026.
I thank Zhelun Chen and Kaiyuan Gu for reading the manuscript and suggesting
revisions.
I also thank Zhelun Chen, Yuxing Ye, Ruida Di, Kaiyuan Gu, and Jiawei Yu for participating
in a seminar on this problem and for their comments.
Finally, I thank Professor Xinyi Yuan for his support during the preparation
of this paper and for supporting my visit to Peking University.

\section*{Use of generative AI}
\addcontentsline{toc}{section}{Use of generative AI}
Generative AI was used in the preparation and revision of this manuscript.
An account of this assistance and selected development records are available
in the supplementary repository:
\begin{quote}
\small\url{https://github.com/DongKaicheng/one-complex-metric-ai-development}.
\end{quote}
The \href{https://github.com/DongKaicheng/one-complex-metric-ai-development/releases/tag/v2026-11-06}{release of 6 November 2026}
includes the
\href{https://github.com/DongKaicheng/one-complex-metric-ai-development/releases/download/v2026-11-06/AI_Assisted_Development_Report_EN.pdf}{English AI-assisted development report (PDF)}
and the
\href{https://github.com/DongKaicheng/one-complex-metric-ai-development/releases/download/v2026-11-06/One_Complex_Metric_AI_Work_Process_GitHub_20261106.zip}{selected mathematical work-process archive (ZIP)}.

\section{Adelic line bundles and complex metrics}
\label{sec:preliminaries}

We use the conventions of \cite{YZ} and additive notation for tensor products of line bundles.

\subsection{Compactified line bundles on quasi-projective varieties}
\label{sec:geometric-conventions}

Let $k$ be a field of characteristic zero, and let $U/k$ be a
smooth geometrically integral quasi-projective variety.
Fix a smooth projective compactification $U\subset X_0$ and
an effective integral Cartier divisor $D_0$ with
$\Supp D_0=X_0\setminus U$.
When $U$ is projective, we take $X_0=U$ and $D_0=0$.

\begin{definition}[Projective models and modifications]\label{def:model}
A projective model of $U$ is a normal projective variety $X$
equipped with a specified dense open immersion $U\subset X$.
We work with projective models dominating $X_0$.
A \emph{projective modification} $p:X\to X_0$ is a projective
birational morphism that restricts to the identity on $U$.
For such a modification, we write
\begin{equation}\label{eq:full-boundary}
D_X=p^*D_0,
\end{equation}
where the pullback retains all Cartier multiplicities.

For a fixed $\QQ$-line bundle $L_U$, a model presentation
is a triple $(X,P,\alpha_P)$, where $P$ is a $\QQ$-line
bundle on $X$ and
$\alpha_P:L_U\xrightarrow{\sim}P|_U$
is a specified isomorphism.

For two model presentations $(X,P,\alpha_P)$ and
$(X',Q,\alpha_Q)$ of $L_U$, choose a common projective
model $U\subset W$, with modifications $r_P:W\to X$ and
$r_Q:W\to X'$.
Their difference is represented by the $\QQ$-Cartier divisor
\begin{equation}\label{eq:comparison-divisor}
\Delta(Q,P)
=
\divisor\bigl(
r_Q^*\alpha_Q\circ(r_P^*\alpha_P)^{-1}
\bigr),
\end{equation}
which is supported on $W\setminus U$.
\end{definition}

\begin{definition}[Boundary norm]\label{def:boundary-norm}
Suppose $D_0\ne0$.
For two model presentations $P$ and $Q$ as above,
let $U\subset W$ be a common projective model, with
induced modification $\rho:W\to X_0$.
Their boundary distance is
\begin{equation}\label{eq:boundary-norm}
\bnorm{Q-P}
=
\inf\bigl\{
c\in\QQ_{>0}:
-c\rho^*D_0
\leq\Delta(Q,P)
\leq c\rho^*D_0
\bigr\}.
\end{equation}
The inequalities are inequalities of $\QQ$-Cartier
divisors. This norm is independent of the chosen common
model; see \cite[Section~2.4]{YZ}.
\end{definition}

\begin{definition}[Compactified line bundles]
\label{def:completion}
A compactified $\QQ$-line bundle $\LL$ on $U/k$
with underlying line bundle $L_U$ is represented by a
boundary-Cauchy sequence of model presentations of $L_U$.
Two sequences represent the same compactified line bundle
if there is an isomorphism between their underlying line
bundles over $U$ for which the boundary norms of their
differences tend to zero.

We denote the group of compactified $\QQ$-line bundles
by $\AdPic(U/k)_{\QQ}$ and identify a model presentation
$(X,P,\alpha_P)$ with its constant class, denoting
both by $P$. When $U$ is projective, we set
$\AdPic(U/k)_{\QQ}=\Pic(U)_{\QQ}$, so that $\LL=L_U$.
\end{definition}

Thus $\AdPic(U/k)_{\QQ}$ is the group of adelic line bundles
in the geometric setting of Yuan--Zhang
\cite[Sections~1.5 and~2.5.1]{YZ}.

When $D_0\ne0$, the boundary norm extends to completed differences
with a specified identification over $U$.
We use the order defined by the closure of the cone of
effective boundary divisors. Thus
\begin{equation}\label{eq:completed-order}
\bnorm{v}\leq c
\quad\Longleftrightarrow\quad
-cD_0\leq v\leq cD_0,
\end{equation}
where $v$ is a completed boundary difference and
$c\in\QQ_{\geq0}$.
Here $D_0$ denotes its constant class with the canonical
trivialization over $U$.

\begin{definition}[Nefness]\label{def:nefness}
Following \cite[Definition~2.5.2]{YZ}, a compactified
line bundle is \emph{strongly nef} if it is a boundary
limit of nef projective model line bundles.

A compactified line bundle $\LL$ is \emph{nef} if there
exists a fixed strongly nef compactified line bundle
$\mathcal A$ such that $a\LL+\mathcal A$ is strongly
nef for every integer $a\geq1$.

A compactified line bundle is \emph{integrable} if it is
a difference of two strongly nef compactified line bundles.
\end{definition}

\subsection{Complex metrics and multiplier ideals}

In this subsection the varieties are defined over $\CC$.
We use the standard terminology for singular Hermitian
metrics and plurisubharmonic functions, as in
\cite{DemCADG,DemBook}.

\begin{definition}[Metric and curvature conventions]
\label{def:weights}
Let $h$ be a possibly singular Hermitian metric with
locally integrable weights on a holomorphic line bundle $L$.
For a nonvanishing local holomorphic frame $e$,
we define its local weight by
$\varphi_{h,e}=-\log\|e\|_h$.

We use the convention
$dd^c=\frac{i}{\pi}\partial\bar\partial$
so that the first Chern current of $(L,h)$ is
\[
c_1(L,h)=dd^c\varphi_{h,e}.
\]
The metric $h$ is called \emph{semipositive} if its first
Chern current $c_1(L,h)$ is a positive $(1,1)$-current;
we write $c_1(L,h)\geq0$.
For singular semipositive metrics, local weights are understood
in their plurisubharmonic (psh) representatives.

For $\QQ$-line bundles, these conventions are
defined by passing to an integral tensor multiple
and dividing by that multiple.
\end{definition}

\begin{definition}[Boundary Green function]
\label{def:boundary-green}
Suppose $D_0\ne0$. Choose a smooth Hermitian metric
on $\OO(D_0)$ and rescale it by a positive constant
so that
\begin{equation}\label{eq:boundary-gauge}
g_0=-\log\|s_{D_0}\|\geq1
\qquad\text{on }U,
\end{equation}
where $s_{D_0}$ is the canonical section of
$\OO(D_0)$.

On a projective model $U\subset X$ equipped with a
modification $p:X\to X_0$, let $g_X=p^*g_0$. Locally, choose a regular frame
$e_{D_X}$ of $\OO(D_X)$ and write
$s_{D_X}=f_{D_X}e_{D_X}$. Then
\begin{equation}\label{eq:local-gauge}
g_X=-\log|f_{D_X}|+\gamma_X,
\end{equation}
where $\gamma_X$ is smooth.
\end{definition}

\begin{definition}[Compactified metrics]\label{def:compactified-metric}
Let $\LL$ be a compactified $\QQ$-line bundle on a smooth
complex quasi-projective variety $U$, with underlying line
bundle $L_U$.
A \emph{compactified metric} on $\LL$ is a continuous Hermitian
metric $h$ on $L_U$ admitting continuous metrized model
presentations satisfying \eqref{eq:simultaneous-approximation},
with respect to the fixed boundary data $X_0,D_0,g_0$.
The algebraic and metric comparisons must use the same
specified isomorphisms over $U$; no positivity is required
of the approximating metrics $h_t$.
If $U$ is projective, every continuous Hermitian metric on
$L_U$ is a compactified metric.
\end{definition}

The condition is independent of the boundary data. On a common
modification, two effective boundary divisors with full boundary
support bound one another up to positive rational constants,
as do their normalized Green functions. Shrinking the tolerance
therefore transfers \eqref{eq:simultaneous-approximation} from
one choice to the other. Pullback to a resolution that is an
isomorphism over $U$ preserves both estimates, so the approximating
modifications may be taken to have smooth sources.

\begin{definition}[Multiplier ideal sheaf]
\label{def:multiplier-ideal}
Let $w$ denote the local psh weights of a singular Hermitian
metric on a line bundle over a smooth complex variety $X$.
For a positive integer $m$, define
\begin{equation}\label{eq:multiplier-definition}
\JJ(mw)_x
=
\left\{
a\in\OO_{X^{\an},x}:
|a|^2e^{-2mw}
\text{ is locally integrable near }x
\right\}.
\end{equation}
Here integrability is understood with respect to Lebesgue
measure in local holomorphic coordinates, or equivalently,
any smooth positive volume form.

The local ideals glue to a coherent analytic ideal sheaf
$\JJ(mw)\subseteq\OO_{X^{\an}}$, called the
\emph{multiplier ideal sheaf}
\cite[Proposition~5.7]{DemBook}.
The local definitions agree on overlaps because the weights
differ by locally bounded functions.

\end{definition}

\subsection{Arithmetic adelic line bundles}

\begin{definition}[Arithmetic adelic line bundles]
\label{def:arithmetic-adelic}
Let $K$ be a number field and let $U/K$ be as above. We write
$\AdPic(U/\ZZ)_{\QQ}$ for the arithmetic adelic
$\QQ$-line bundles of
\cite[Section~2.5]{YZ}.
Such a line bundle is represented by a Cauchy sequence
of Hermitian line bundles on projective arithmetic models.

More precisely, relative to a strictly effective arithmetic
boundary $\overline{\mathscr D}_0$, the comparisons in
a cofinal model presentation satisfy
\begin{equation}\label{eq:arithmetic-cauchy}
-\epsilon_i\overline{\mathscr D}_0
\leq
\widehat{\divisor}(\alpha_j\alpha_i^{-1})
\leq
\epsilon_i\overline{\mathscr D}_0
\qquad (j\geq i),
\end{equation}
for positive rational numbers $\epsilon_i\to0$.
These are arithmetic divisor inequalities, including
their Green-function components.
\end{definition}

\begin{definition}[Geometric part and complex metric]
\label{def:geometric-part}
For $\overline{\LL}\in\AdPic(U/\ZZ)_{\QQ}$,
its \emph{geometric part}
$\LL^{\geom}\in\AdPic(U/K)_{\QQ}$
is obtained by taking the generic fibers of a model
presentation, retaining the specified isomorphisms over $U$.

For every embedding $\sigma:K\hookrightarrow\CC$,
the same arithmetic presentation determines a continuous
Hermitian metric $h_\sigma$ on the underlying line bundle
over $U_\sigma$. These constructions are independent of
the presentation; see
\cite[Section~2.5.5 and Proposition~3.4.1]{YZ}.
\end{definition}

Arithmetic nefness is understood in the sense of
\cite[Section~2.5]{YZ}. In particular, it implies
nefness of the geometric part. The converse is not asserted:
the geometric nefness considered in this paper does not
include the arithmetic degree conditions.

\subsection{The Gross--Schoen line bundles}

We recall the definitions of Gross--Schoen adelic line
bundles given by Yuan \cite[Section~3.3.5]{Yuan} and
Gao--Zhang \cite[Section~4.1]{GZ}.

\begin{definition}[Yuan's Gross--Schoen line bundle]
\label{def:yuan-GS}
Let $S$ be a smooth geometrically integral quasi-projective
variety over a number field, and let $\pi:\mathcal C\to S$
be a smooth projective family of curves of genus $g\geq2$.
Write $\overline\omega_a$ and
$\overline{\OO}(\Delta)_a$ for the canonical admissible
adelic extensions of the relative dualizing and
diagonal line bundles. Set
\begin{equation}\label{eq:phi-line}
\overline{\OO}(\Phi_S)
=
\langle\overline\omega_a,\overline\omega_a\rangle_\pi
-
\langle
\overline{\OO}(\Delta)_a,
\overline{\OO}(\Delta)_a,
\overline{\OO}(\Delta)_a
\rangle_{\pi^{(2)}},
\end{equation}
where $\pi^{(2)}:\mathcal C\times_S\mathcal C\to S$.
Its underlying line bundle has a canonical trivialization
induced by adjunction.

In this notation, Yuan's Gross--Schoen adelic line bundle is
\begin{equation}\label{eq:yuan-GS}
\overline G^{\,Y}_{\mathcal C/S}
=
\frac{2g+1}{2g-2}
\langle\overline\omega_a,\overline\omega_a\rangle_\pi
-
\overline{\OO}(\Phi_S).
\end{equation}
Here $\overline{\OO}(\Phi_S)$ is the adelic line bundle
associated with the relative $\varphi$-invariant in
\cite[Section~3.3.2]{Yuan}; the Gross--Schoen combination
\eqref{eq:yuan-GS} is introduced in
\cite[Section~3.3.5]{Yuan}.
\end{definition}

\begin{definition}[Gao--Zhang's Gross--Schoen line bundle]
\label{def:GZ-GS}
Let $S$ be a smooth quasi-projective variety over a
number field, and let $\pi:\mathcal C\to S$ be a smooth
projective family of curves of genus $g\geq3$.
Suppose that the family is equipped with a relative
degree-one class $\xi\in\operatorname{Pic}^1(\mathcal C/S)$
satisfying
\[
(2g-2)\xi=\omega_{\mathcal C/S}.
\]

Following \cite[Section~4.1]{GZ}, let
$J=\operatorname{Jac}(\mathcal C/S)$ and
$i_\xi:\mathcal C\to J$ be the Abel--Jacobi morphism.
Let $\overline{\mathcal P}^{\Delta}$ be the canonical
adelic extension of the diagonal Poincar\'e line bundle
on $J$, characterized by
\[
[2]^*\overline{\mathcal P}^{\Delta}
=4\overline{\mathcal P}^{\Delta}.
\]
Writing $m,p_1,p_2:J\times_SJ\to J$ for the addition
morphism and the two projections, respectively, define
\begin{equation}\label{eq:poincare-adelic}
2\overline{\mathcal P}
=
m^*\overline{\mathcal P}^{\Delta}
-p_1^*\overline{\mathcal P}^{\Delta}
-p_2^*\overline{\mathcal P}^{\Delta},
\qquad
\overline{\mathcal Q}
=
(i_\xi,i_\xi)^*\overline{\mathcal P}.
\end{equation}

Gao--Zhang's Gross--Schoen adelic line bundle is
\begin{equation}\label{eq:GZ-GS}
\overline G_g
=
\langle
\overline{\mathcal Q},
\overline{\mathcal Q},
\overline{\mathcal Q}
\rangle_{\pi^{(2)}},
\end{equation}
where $\pi^{(2)}:\mathcal C\times_S\mathcal C\to S$.
We denote its geometric part by $G_g^{\geom}$ and its
underlying line bundle on $S$ by $G_{g,S}$.

This is the pullback of Gao--Zhang's construction on the
moduli stack of pairs $(C,\xi)$
\cite[Section~4.1]{GZ}.
\end{definition}

The height and curvature of this metric are recalled in
\Cref{prop:GS-metric}.

\section{Proof of Theorem~\ref{thm:main-nefness}}\label{sec:main-proof}
We prove \Cref{thm:main-nefness} together with the semiample
approximation statement of \Cref{rem:semiample-approximation}.
Throughout this section, the varieties are defined over $\CC$
and the conventions of \Cref{sec:geometric-conventions} apply.
Since the compactified metric is continuous and has semipositive
curvature, its local weights are continuous psh functions.

\subsection{Extension across the boundary}\label{sec:corrections}

\begin{theorem}[Removable singularities for plurisubharmonic functions]
\label{thm:psh-extension}
Let $V$ be a connected complex manifold and
$Z\subset V$ a proper closed analytic subset.
A psh function on $V\setminus Z$ that is locally
bounded above near $Z$ has a unique psh extension,
obtained by upper limiting values.
The extension is locally integrable unless identically
$-\infty$.

Moreover, adding a pluriharmonic function on $V$ commutes
with this extension.
This is the removable-singularity theorem of
\cite[Chapter~I, Theorem~5.24]{DemCADG}.
\end{theorem}

\Needspace{7\baselineskip}
\begin{proposition}\label{prop:small-extension}
Assume that $D_0\ne0$ and that $(\LL,h)$ satisfies the
hypotheses of \Cref{thm:main-nefness}.
For a rational $t>0$, choose a projective modification
$p_t:X_t\to X_0$ with $X_t$ smooth, together with
$P_t,h_t,\alpha_t$ satisfying
\eqref{eq:simultaneous-approximation}, and put
\[
D_t=p_t^*D_0,
\qquad
B_t=P_t+t D_t.
\]
Identify $B_t|_U$ with $L_U$ using $\alpha_t$
and the canonical section of $t D_t$.

Then $h$ extends to a singular semipositive metric on $B_t$.
Moreover,
\begin{equation}\label{eq:corrected-completion}
0\leq B_t-\LL\leq2t D_0
\end{equation}
in the completed boundary order.

On every sufficiently small regular chart of $X_t$,
the local weight $w_t$ of the extended metric satisfies
\begin{equation}\label{eq:small-singularity}
2t\log|f_{D_t}|-C
\leq w_t\leq C,
\end{equation}
where $f_{D_t}$ is a local equation of $D_t$
and $C$ is a constant depending on the chart.
\end{proposition}

\begin{proof}
Locally, choose regular frames of $P_t$ and
$\OO(D_t)$, and let $f_{D_t}$ be the local equation
of $D_t$. Write $\varphi$ and $\beta_t$ for the
weights of $h$ and $h_t$, respectively, in the
corresponding frames. The metric induced by $h$ on
$B_t|_U$ has weight
\begin{equation}\label{eq:transported-weight}
w_t=\varphi+t\log|f_{D_t}|.
\end{equation}

Since $g_{X_t}=-\log|f_{D_t}|+\gamma_t$,
condition \eqref{eq:simultaneous-approximation} gives
\[
\beta_t+2t\log|f_{D_t}|-t\gamma_t
\leq w_t
\leq\beta_t+t\gamma_t.
\]
As $\beta_t$ and $\gamma_t$ are locally bounded,
\eqref{eq:small-singularity} follows.
In particular, $w_t$ is locally bounded above near
the boundary. By \Cref{thm:psh-extension}, it extends
uniquely as a psh function across $D_t$.
These local extensions glue by uniqueness, giving
the required singular semipositive metric on $B_t$.

Finally, \eqref{eq:simultaneous-approximation} implies
$-t D_0\leq P_t-\LL\leq t D_0$
in the completed boundary order.
Adding $t D_0$ yields
\eqref{eq:corrected-completion}.
\end{proof}

\subsection{Ample perturbations on a projective model}
\label{sec:ample-perturbations}

For a projective morphism $p:X\to X_0$, twisting a $p$-ample
line bundle by the pullback of a line bundle on $X_0$ preserves
relative ampleness. If $H$ is $p$-ample and $A$ is ample on
$X_0$, then $H+np^*A$ is ample for all sufficiently large
integers $n$ \cite[Tags~01VJ and~0892]{Stacks}.

\begin{lemma}[Ample perturbation]
\label{lem:ample-adjustment}
Let $p:X\to X_0$ be a projective modification with $X$ smooth,
and let $A$ be an ample $\QQ$-line bundle on $X_0$.
For every rational $\lambda>0$, there exists an ample
$\QQ$-line bundle $A_X$ on $X$, equipped with an
identification $A_X|_U\simeq A|_U$, such that
\begin{equation}\label{eq:ample-adjustment}
-\lambda D_X\leq A_X-p^*A\leq\lambda D_X.
\end{equation}
Here the inequalities are inequalities of boundary
$\QQ$-Cartier divisors determined by the specified
identification over $U$.
\end{lemma}

\begin{proof}
Choose an integral $p$-ample line bundle $H$ on $X$.
Since $X_0$ is smooth, the line bundle $H|_U$ extends
to a line bundle $E$ on $X_0$
\cite[Tag~0BD9]{Stacks}.

Then $R=H-p^*E$ is $p$-ample and has a specified
trivialization over $U$. Thus $R$ is represented by
a Cartier divisor supported on $X\setminus U$.
Since $D_X$ has positive multiplicity along every
boundary prime divisor, there exists $c>0$ such that
\[
-cD_X\leq R\leq cD_X.
\]

By relative ampleness, $R+np^*A$ is ample for every
sufficiently large divisible integer $n$.
Choose such an $n$ with $c/n\leq\lambda$, and set
\[
A_X=p^*A+\frac{1}{n}R.
\]
Then $A_X$ is ample, its restriction to $U$ is
identified with $A|_U$, and
\[
-\lambda D_X
\leq \frac{1}{n}R
=A_X-p^*A
\leq\lambda D_X.
\]
\end{proof}

For the projective models $U\subset X_t$ appearing in the proof of
\Cref{thm:main-nefness}, we denote the corresponding
ample perturbations by $A_t$.
\subsection{Multiplier ideals and principalization}\label{sec:ideals}
Fix a rational $t>0$ and the data of \Cref{prop:small-extension}.
For a coherent ideal $\II\subseteq\OO_X$ and a morphism
$\rho:Y\to X$, write $\II\OO_Y=\im(\rho^*\II\to\OO_Y)$
for the image ideal.

\begin{proposition}[Algebraization and principalization of the multiplier ideal]
\label{prop:ideal-principalization}
Choose a positive integer $m$ such that $mB_t$ is integral
and $mt\in\ZZ$, and set $b_m=2mt$.
The multiplier ideal $\JJ(mw_t)$ algebraizes to a
nonzero coherent ideal $\II_{t,m}\subseteq\OO_{X_t}$
satisfying
\begin{equation}\label{eq:ideal-sandwich}
\OO_{X_t}(-b_mD_t)
\subseteq\II_{t,m}\subseteq\OO_{X_t},
\qquad
\II_{t,m}|_U=\OO_U.
\end{equation}
Moreover, there exists a projective modification
$\rho_{t,m}:Y_{t,m}\to X_t$, with $Y_{t,m}$ smooth
and $\rho_{t,m}$ an isomorphism over $U$, such that
\begin{equation}\label{eq:principalized-ideal}
\II_{t,m}\OO_{Y_{t,m}}
=\OO_{Y_{t,m}}(-F_{t,m}),
\end{equation}
where $F_{t,m}$ is an effective Cartier divisor
supported on $Y_{t,m}\setminus U$.
\end{proposition}

\begin{proof}
Let $f$ be a local equation of $D_t$.
By \eqref{eq:small-singularity}, we have
$w_t\geq2t\log|f|-C$.
Since $b_m=2mt$, every local section $f^{b_m}a$
of $\OO(-b_mD_t)$ satisfies
\[
|f^{b_m}a|^2e^{-2mw_t}
\leq e^{2mC}|a|^2,
\]
which is locally integrable. Hence
$\OO(-b_mD_t)^{\an}\subseteq\JJ(mw_t)$.
Moreover, $w_t$ is continuous on $U$, so
$\JJ(mw_t)|_{U^{\an}}=\OO_{U^{\an}}$.

The multiplier ideal is coherent
\cite[Proposition~5.7]{DemBook}.
By projective GAGA
\cite[no.~12, Th\'eor\`emes~1--3]{Serre},
it algebraizes, together with its inclusions, to a
coherent ideal $\II_{t,m}$ satisfying
\eqref{eq:ideal-sandwich}.

Principalization in characteristic zero
\cite[Example~1.19(1)]{BM} gives a projective modification
$\rho_{t,m}$ with smooth source, isomorphic over $U$, on which
the image ideal is invertible. Since $\II_{t,m}|_U=\OO_U$,
its divisor $F_{t,m}$ is effective and supported outside $U$.
\end{proof}
\begin{proposition}[The divisor of a principalized multiplier ideal]\label{prop:ideal-divisor}
Let $X_t,D_t,B_t,w_t,m$ and $\II_{t,m}$ be as in \Cref{prop:ideal-principalization}. Let $\rho:Y\to X_t$ be \emph{any} projective principalization with $Y$ smooth that is an isomorphism over $U$, and let $F$ be the effective Cartier divisor defined by the image ideal
\[
 \II_{t,m}\OO_Y=\OO_Y(-F).
\]
Then
\begin{equation}\label{eq:exceptional-bound}
 0\leq\frac{F}{m}\leq2t\rho^*D_t
\end{equation}
as $\QQ$-Cartier divisors on $Y$. 
\end{proposition}
\begin{proof}
Fix $y\in Y$ and put $x=\rho(y)$. Choose a local equation $f$ of $D_t$ at $x$, and local generators $a_1,\ldots,a_s$ of $\II_{t,m,x}$. The first inclusion in \eqref{eq:ideal-sandwich} means that there are $c_i\in\OO_{X_t,x}$ with
\[
 f^{b_m}=\sum_{i=1}^s c_i a_i,\qquad b_m=2mt.
\]
Apply the local ring homomorphism $\rho^\#: \OO_{X_t,x}\to\OO_{Y,y}$. Its image gives
\[
 (\rho^\#f)^{b_m}\in \II_{t,m}\OO_{Y,y}.
\]
This uses the image ideal and requires no flatness assumption on $\rho$.

Choose a nonzero local generator $g$ of $\OO_Y(-F)$.
The membership above gives
\begin{equation}\label{eq:local-divisibility}
 (\rho^\#f)^{b_m}=g h_y
\end{equation}
for a regular function $h_y$ on a sufficiently small neighborhood.
Since $\rho$ is dominant, $h_y\ne0$, and hence
\[
 b_m\rho^*D_t-F=\divisor(h_y)\geq0.
\]
Together with $F\geq0$, this proves $0\leq F\leq b_m\rho^*D_t$.
Dividing by $m$ gives \eqref{eq:exceptional-bound}.
\end{proof}

\subsection{Global generation}\label{sec:approximation}

\begin{theorem}[Nadel vanishing and Castelnuovo--Mumford regularity]\label{thm:nadel-regularity}
Let $X$ be a smooth projective complex variety.
\begin{enumerate}
\item If an integral line bundle $T$ carries a singular Hermitian metric $h_T$ with curvature bounded below by a smooth positive form, then
\[
 H^i\bigl(X,K_X\otimes T\otimes\JJ(h_T)\bigr)=0
 \quad\text{for every }i>0.
\]
The ideal is defined from the local weights of $h_T$ as in \Cref{def:multiplier-ideal}.
\item If $H_0$ is very ample and a coherent sheaf $\mathscr S$ satisfies $H^i(X,\mathscr S\otimes H_0^{-i})=0$ for every $i>0$, then $\mathscr S$ is globally generated.
\end{enumerate}
The first assertion is analytic Nadel vanishing
\cite[Theorem~5.9]{BBGH}; see also \cite{Nadel}.
The second is Castelnuovo--Mumford regularity
\cite[Tag~08A8]{Stacks}, applied to the embedding
defined by $H_0$.
Projective GAGA identifies analytic and algebraic
coherent cohomology.
\end{theorem}

\begin{lemma}[Generation after an ample perturbation]\label{lem:exact-generation}
Fix rational $t,\varepsilon>0$, the line bundle $B_t$ with its semipositive metric from \Cref{prop:small-extension}, and an ample line bundle $A_t$ as in \Cref{lem:ample-adjustment}. Give $A_t$ a smooth positive metric and set
\[
 E_{t,\varepsilon}=B_t+\varepsilon A_t.
\]
For all sufficiently large integers $m$ divisible enough that $mE_{t,\varepsilon}$ and $mB_t$ are integral and $mt\in\ZZ$, the coherent sheaf
\begin{equation}\label{eq:section-sheaf}
 \mathscr S_{t,\varepsilon,m}
  =\OO_{X_t}(mE_{t,\varepsilon})\otimes\II_{t,m}
\end{equation}
is globally generated. The threshold depends on the fixed data.
\end{lemma}
\begin{proof}
The product metric on $E_{t,\varepsilon}$ has curvature at
least $\varepsilon\omega_{A_t}>0$. Put $n=\dim X_t$.
Choose an integral very ample line bundle $H_0$ with a smooth
positive metric, and a smooth metric on $K_{X_t}$.
For $1\leq i\leq n$, consider
\[
T_{m,i}=mE_{t,\varepsilon}-K_{X_t}-iH_0.
\]
On the fixed compact variety $X_t$, there is a constant
$C$ such that
\[
c_1(K_{X_t})+n\omega_{H_0}\leq C\omega_{A_t}.
\]
The induced singular metric on $T_{m,i}$ therefore satisfies
\[
c_1(T_{m,i})\geq(m\varepsilon-C)\omega_{A_t}>0
\]
for $m\varepsilon>C$. Its weights differ from $mw_t$
by smooth functions, so its multiplier ideal is
$\II_{t,m}^{\an}$.

Nadel vanishing and GAGA give
\[
H^i\bigl(X_t,
\OO(mE_{t,\varepsilon}-iH_0)\otimes\II_{t,m}\bigr)=0
\qquad (1\leq i\leq n).
\]
For $i>n$, cohomology vanishes by dimension.
Castelnuovo--Mumford regularity now implies that
$\mathscr S_{t,\varepsilon,m}$ is globally generated.
\end{proof}

\subsection{Completion of the proof}
\label{sec:completion-proof}
\phantomsection\label{proof:main-nefness}

\begin{proof}[Proof of \Cref{thm:main-nefness,rem:semiample-approximation}]
Fix $X_0,D_0$ and an ample $\QQ$-line bundle $A$ on $X_0$.
Let $h$ be a compactified metric on $\LL$ with
$c_1(L_U,h)\geq0$.
Assume first that $D_0\ne0$, and fix a rational
$\varepsilon>0$.

\paragraph{Construction of the model line bundle.}
For a rational $t>0$, choose a projective modification
$p_t:X_t\to X_0$ with $X_t$ smooth and approximation
data satisfying \eqref{eq:simultaneous-approximation}.
By \Cref{prop:small-extension}, the line bundle $B_t$
admits a singular semipositive metric and satisfies
\[
0\leq B_t-\LL\leq2t D_0.
\]

Choose a rational $0<\lambda\leq t/(1+\varepsilon)$.
By \Cref{lem:ample-adjustment}, there exists an ample
$\QQ$-line bundle $A_t$ on $X_t$ such that
\[
-\lambda D_t
\leq A_t-p_t^*A
\leq\lambda D_t.
\]
Equip $A_t$ with a smooth positive metric.

With these choices fixed, take a sufficiently large divisible
integer $m$ satisfying the conditions of
\Cref{lem:exact-generation}.
By \Cref{prop:ideal-principalization}, there is a projective
modification $\rho=\rho_{t,m}:Y_{t,m}\to X_t$
with $Y_{t,m}$ smooth, such that
$\II_{t,m}\OO_{Y_{t,m}}
=\OO_{Y_{t,m}}(-F_{t,m})$.
For simplicity, write $F=F_{t,m}$ and
$q=p_t\circ\rho$.

Define
\begin{equation}\label{eq:candidate-model}
N_{t,\varepsilon,m}
=\rho^*(B_t+\varepsilon A_t)-\frac{F}{m}.
\end{equation}

\paragraph{Global generation and restriction to $U$.}
Set
\[
 V_m=H^0(X_t,\mathscr S_{t,\varepsilon,m}).
\]
By \Cref{lem:exact-generation}, the evaluation map
\[
 V_m\otimes\OO_{X_t}\twoheadrightarrow\mathscr S_{t,\varepsilon,m}
\]
is surjective. Pulling back to $Y_{t,m}$ and passing to the image
ideal $\II_{t,m}\OO_{Y_{t,m}}=\OO_{Y_{t,m}}(-F)$ gives
\[
V_m\otimes\OO_{Y_{t,m}}
\twoheadrightarrow
\OO_{Y_{t,m}}\bigl(m\rho^*(B_t+\varepsilon A_t)-F\bigr)
=\OO_{Y_{t,m}}(mN_{t,\varepsilon,m}).
\]
Hence $N_{t,\varepsilon,m}$ is semiample.

Since $\II_{t,m}|_U=\OO_U$ and the modifications
are isomorphisms over $U$, the specified identifications
for $B_t$ and $A_t$ induce
\[
N_{t,\varepsilon,m}|_U
\simeq L_U+\varepsilon A|_U.
\]
All subsequent boundary comparisons use this identification.

\paragraph{Boundary estimate.}
Put $D_Y=q^*D_0=\rho^*D_t$.
By \Cref{prop:ideal-divisor},
\[
0\leq\frac{F}{m}\leq2t D_Y.
\]
Together with \Cref{prop:small-extension}, this gives
\[
-2t D_0
\leq
(\rho^*B_t-\LL)-\frac{F}{m}
\leq2t D_0
\]
in the completed boundary order.

Using \eqref{eq:candidate-model}, we have
\begin{align*}
N_{t,\varepsilon,m}-(\LL+\varepsilon A)
={}&(\rho^*B_t-\LL)-\frac{F}{m}\\
&+\varepsilon\rho^*(A_t-p_t^*A).
\end{align*}
The ample adjustment satisfies
\[
-\varepsilon\lambda D_Y
\leq\varepsilon\rho^*(A_t-p_t^*A)
\leq\varepsilon\lambda D_Y.
\]
Consequently,
\begin{equation}\label{eq:final-boundary-error}
\bnorm{N_{t,\varepsilon,m}-(\LL+\varepsilon A)}
\leq2t+\varepsilon\lambda
\leq3t.
\end{equation}
This estimate holds for every sufficiently large admissible $m$.

\paragraph{Convergence and nefness.}
Given $\delta>0$, choose a rational
$0<t<\delta/3$ and construct
$N_{t,\varepsilon,m}$ as above.
Then \eqref{eq:final-boundary-error} proves
\eqref{eq:main-approximation}.

To obtain a convergent presentation, choose
$t_j\downarrow0$ and an admissible $m_j$ for each
$t_j$. By the triangle inequality,
\[
\bnorm{
N_{t_j,\varepsilon,m_j}
-
N_{t_k,\varepsilon,m_k}
}
\leq3t_j+3t_k.
\]
Thus the semiample model line bundles
$N_{t_j,\varepsilon,m_j}$ form a boundary-Cauchy
sequence converging to $\LL+\varepsilon A$.
It follows that $\LL+\varepsilon A$ is strongly nef
for every rational $\varepsilon>0$.

For any integer $a\geq1$, taking $\varepsilon=1/a$
and multiplying the resulting presentation by $a$
shows that $a\LL+A$ is strongly nef.
Since $A$ is a fixed strongly nef line bundle,
independent of $a$, \Cref{def:nefness} implies
that $\LL$ is nef.

If $U$ is projective, then $D_0=0$ and $\LL=L_U$.
For every rational $\varepsilon>0$, the line bundle
$L_U+\varepsilon A$ has a continuous metric with strictly
positive curvature. The Nadel--regularity argument of
\Cref{lem:exact-generation}, with trivial multiplier ideals,
shows that it is semiample. Hence $\LL$ is nef in this case
as well.
\end{proof}

\section{Approximation over a number field}
\label{sec:arithmetic}

\subsection{Arithmetic adelic presentations}

\begin{proposition}[Compatibility of geometric and complex realizations]
\label{prop:adelic-approximation}
Let $K$ be a number field, and let $U/K$ be smooth,
geometrically integral, and quasi-projective.
For any arithmetic adelic $\QQ$-line bundle
$\overline{\LL}\in\AdPic(U/\ZZ)_{\QQ}$ and any embedding
$\sigma:K\hookrightarrow\CC$, its prescribed continuous metric $h_\sigma$ is a compactified
metric on the geometric part $(\LL^{\geom})_\sigma$.
No positivity or integrability assumption on
$\overline{\LL}$ is required.
\end{proposition}

\begin{proof}
This follows from the geometric-part construction
\cite[Section~2.5.5]{YZ} and the convergence of model
metrics in the proof of
\cite[Proposition~3.4.1, Section~3.4.2]{YZ},
together with the independence of the boundary topology
\cite[Section~2.4]{YZ}.
\end{proof}

\subsection{Specialization and descent}

\begin{proposition}[Approximation over the field of definition]
\label{prop:finite-transfer}
Let $k\subset\CC$ be a field of characteristic zero, let
$X/k$ be smooth, projective, and geometrically integral, and
let $U\subset X$ be a nonempty open subset.
Let $D$ be an effective integral Cartier divisor with
$\Supp D=X\setminus U$, and let $Q\in\Pic(X)_{\QQ}$.
The case $U=X$ is allowed, with $D=0$.

Suppose that, for some rational $\delta>0$, there are a
projective modification $p_{\CC}:Y_{\CC}\to X_{\CC}$,
with $Y_{\CC}$ normal and $p_{\CC}$ an isomorphism over
$U_{\CC}$, and a semiample $\QQ$-line bundle
$N_{\CC}$ on $Y_{\CC}$, with a specified identification
$N_{\CC}|_{U_{\CC}}\simeq Q|_{U_{\CC}}$, such that
\begin{equation}\label{eq:transfer-assumption}
-\delta p_{\CC}^*D
\leq N_{\CC}-p_{\CC}^*Q
\leq\delta p_{\CC}^*D.
\end{equation}
Then there exist a projective modification $p:Z\to X$ over $k$,
with $Z$ smooth and geometrically integral and $p$ an isomorphism
over $U$, and a semiample
$\QQ$-line bundle $N_0$ on $Z$, with a specified identification
$N_0|_U\simeq Q|_U$, such that
\begin{equation}\label{eq:transfer-conclusion}
-\delta p^*D\leq N_0-p^*Q\leq\delta p^*D.
\end{equation}
The inequalities are inequalities of comparison divisors
determined by the specified identifications over $U$.
\end{proposition}

\begin{proof}
Choose a positive integer $b$ clearing all denominators
and realizing the open isomorphism, such that
$H_{\CC}=bN_{\CC}$ is globally generated.
The inequalities \eqref{eq:transfer-assumption} give regular maps
\[
p_{\CC}^*\OO_{X_{\CC}}(bQ_{\CC}-b\delta D_{\CC})
\longrightarrow H_{\CC}
\longrightarrow
p_{\CC}^*\OO_{X_{\CC}}(bQ_{\CC}+b\delta D_{\CC}),
\]
whose restrictions to $U_{\CC}$ are determined by the
specified identification and the canonical boundary
trivialization.

By spreading out
\cite[Tags~01ZM, 01ZR, and~0B8W]{Stacks}, the projective
modification, line bundle, generating sections, specified
isomorphisms over $U$, and comparison maps
are defined over a finitely generated integral
$k$-subalgebra of $\CC$.
After localization, specialize at a closed point with
residue field $\ell$, finite over $k$.
Replace the fiber by the normalization of the reduced
closure of $U_\ell$; this gives a projective modification
$Y_\ell\to X_\ell$, unchanged over $U_\ell$
\cite[Tags~0BXR and~0GK4]{Stacks}.
The generating sections and comparison maps pull back,
and the latter remain nonzero on $U_\ell$.
Thus the resulting globally generated line bundle
$H_\ell$ satisfies the same bounds, with tolerance
$b\delta$.

Choose a finite Galois extension $F/k$ containing $\ell$,
with group $G$ and degree $r$.
Let $W_F$ be the normalization of the closure of the
diagonal $U_F$ in
\[
\prod_{\tau\in G,\,X_F}{}^\tau Y_F,
\]
and let $r_\tau:W_F\to{}^\tau Y_F$ be the projections.
Permutation of the factors gives a Galois descent datum.
The tensor product of the conjugate pullbacks of an
auxiliary very ample line bundle on $Y_F$ gives a
compatible ample line bundle on $W_F$.
By Galois descent and descent of ampleness
\cite[Tags~0CCJ, 0CDQ, and~0D2P]{Stacks}, $W_F$ descends
to a projective modification $p:Z\to X$ with $Z$ normal,
isomorphic over $U$.
The same descent applies to
\[
\mathscr B_F
=\bigotimes_{\tau\in G}r_\tau^*({}^\tau H_F),
\]
giving a line bundle $\mathscr B$ on $Z$ with a specified
identification $\mathscr B|_U\simeq brQ|_U$.
It is globally generated by
\cite[Tag~0B57]{Stacks}.

Tensoring and descending the comparison maps gives
\[
-br\delta p^*D
\leq\mathscr B-brp^*Q
\leq br\delta p^*D.
\]
Hence
\begin{equation}\label{eq:normalized-line}
N_0=\frac{1}{br}\mathscr B
\end{equation}
is semiample and satisfies \eqref{eq:transfer-conclusion}
with the required identification over $U$.
Finally, resolve $Z$ without changing $U$
\cite[Theorem~6.1]{BM} and pull back $N_0$.
The resolved projective model $U\subset Z$ is geometrically
integral because it is birational to $X$.
\end{proof}

\subsection{Proof of the number-field theorem}
\label{sec:proof-number-field}
\phantomsection\label{proof:number-field}

\begin{proof}[Proof of \Cref{thm:number-field}]
Choose a smooth projective $K$-compactification $U\subset X_0$,
an effective integral Cartier divisor $D_0$ with support
$X_0\setminus U$, and an ample $\QQ$-line bundle $A$ on $X_0$.
Assume first that $D_0\ne0$.
Fix rational $\varepsilon>0$ and real $\delta>0$, and
choose rational $t>0$ with $3t<\delta$.
Choose a model presentation $Q$ of $L_U$ over $K$ with
$\bnorm{Q-\LL}<t$.
After resolution, we may suppose that $Q$ is defined
on a smooth projective model $U\subset X$, with a
modification $p:X\to X_0$; put $D_X=p^*D_0$.

By \Cref{rem:semiample-approximation}, there is a semiample complex
model line bundle $N_{\CC}$, with the prescribed
restriction to $L_{U,\sigma}+\varepsilon A_\sigma|_{U_\sigma}$,
such that
\[
\bnorm{N_{\CC}-(\LL_\sigma+\varepsilon A_\sigma)}<t.
\]
On a common projective model $U_\sigma\subset W$
dominating $X_\sigma$,
the comparison divisor therefore satisfies
\[
-2tD_W
\leq N_{\CC,W}-(Q+\varepsilon p^*A)_{\sigma,W}
\leq2tD_W.
\]
Here $D_W$ is the full pullback of $D_0$.
Applying \Cref{prop:finite-transfer} with reference
line bundle $Q+\varepsilon p^*A$ and tolerance $2t$
gives a projective modification $r:Y\to X$ over $K$ with
$Y$ smooth and a semiample line bundle $N$, with the required identification
over $U$, such that
\[
-2t r^*D_X
\leq N-r^*(Q+\varepsilon p^*A)
\leq2t r^*D_X.
\]
Consequently,
\begin{equation}\label{eq:three-error}
\begin{aligned}
\bnorm{N-(\LL+\varepsilon A)}
&\leq
\bnorm{N-r^*(Q+\varepsilon p^*A)}
+\bnorm{Q-\LL}\\
&<3t<\delta.
\end{aligned}
\end{equation}

Letting $\delta$ tend to zero gives a boundary-Cauchy
sequence of semiample model line bundles converging
to $\LL+\varepsilon A$, which is therefore strongly nef.
For integers $a,j\geq1$, take $\varepsilon=1/a$ and
$\delta=1/(aj)$, and multiply the approximating line
bundle by $a$. The resulting semiample line bundle
$N_{a,j}$ has the specified restriction $aL_U+A|_U$ and
\begin{equation}\label{eq:one-place-models}
\bnorm{N_{a,j}-(a\LL+A)}<j^{-1}.
\end{equation}
Thus $a\LL+A$ is strongly nef for every $a\geq1$.
Since $A$ is fixed and strongly nef, $\LL$ is nef.

If $U$ is projective, then $\LL=L_U$.
The projective case of \Cref{rem:semiample-approximation}
and descent of global generation \cite[Tag~0B57]{Stacks}
imply that $L_U+\varepsilon A$ is semiample over $K$,
so the same conclusion holds
with zero boundary error.
The assertion for arithmetic adelic line bundles follows
from \Cref{prop:adelic-approximation}.
\end{proof}

\begin{remark}[Semiample approximation over $K$]
\label{rem:number-field-approximation}
Under the hypotheses of \Cref{thm:number-field}, fix a smooth
projective $K$-compactification $U\subset X_0$, an effective
integral Cartier divisor $D_0$ with support $X_0\setminus U$,
and an ample $\QQ$-line bundle $A$ on $X_0$.
For every rational $\varepsilon>0$ and every real $\delta>0$,
the proof gives a projective modification $q:Y\to X_0$ over $K$,
with $Y$ smooth and $q$ an isomorphism over $U$, and a semiample
$\QQ$-line bundle $N$ on $Y$, with a specified isomorphism
\[
 L_U+\varepsilon A|_U\xrightarrow{\sim}N|_U,
\]
such that
\[
 \bnorm{N-(\LL+\varepsilon A)}<\delta.
\]
Thus $\LL+\varepsilon A$ is strongly nef for every rational
$\varepsilon>0$. Taking $\varepsilon=1/a$ and scaling by $a$
shows that $a\LL+A$ is strongly nef for every integer $a\geq1$,
with the same $A$; the approximations satisfy
\eqref{eq:one-place-models}.
If $U$ is projective, one may take $Y=X_0=U$ and
$N=L_U+\varepsilon A$.
\end{remark}

\section{Proof of the Gross--Schoen theorem}\label{sec:GS}
Let $M,\mathcal C,J$, and $\overline G_g$ be as in
\Cref{thm:GS,def:GZ-GS}.

\begin{proposition}[The height and curvature of $\overline G_g$]\label{prop:GS-metric}
The adelic line bundle $\overline G_g$ is integrable. For every algebraic point $s$ of $M$, its normalized point height is
\[
 h_{\overline G_g}(s)
   =\langle\Delta_{\xi_s},\Delta_{\xi_s}\rangle_{\mathrm{BB}}.
\]
At every complex embedding $\sigma:K\hookrightarrow\CC$, its prescribed metric is smooth and satisfies
\begin{equation}\label{eq:GS-curvature}
 c_1(\overline G_{g,\sigma})=\beta_{\GS}\geq0,
\end{equation}
where $\beta_{\GS}$ is the Betti form of the Gross--Schoen normal function, with the Poincar\'e normalization in \Cref{def:GZ-GS}.
\end{proposition}
\begin{proof}
The construction in \cite[Section~4.1]{GZ} gives an integrable adelic Poincar\'e line bundle. Pullback and the relative Deligne pairing preserve integrability by \cite[Section~2.5.5 and Theorem~4.1.3]{YZ}. The height identity is \cite[Theorem~4.1]{GZ}. The smoothness and curvature identity are \cite[Section~4.2 and Theorem~4.2]{GZ}, where the metric is related locally to the metrized biextension associated with the Gross--Schoen normal function. These identities pull back to $M$ and apply to
\eqref{eq:GS-introduction}.
\end{proof}

\phantomsection\label{proof:GS}
\begin{proof}[Proof of \Cref{thm:GS}]
By \Cref{prop:GS-metric}, the prescribed metric of $\overline G_g$
at a complex embedding is smooth and has semipositive first
Chern current. The arithmetic adelic case of
\Cref{thm:number-field} gives the nefness of $G_g^{\geom}$.
\end{proof}

\begin{remark}[Approximation of the Gross--Schoen line bundle]
\label{rem:GS-approximation}
Let $G_{g,M}$ denote the underlying line bundle of
$\overline G_g$ on $M$.
Fix a smooth projective $K$-compactification $M\subset X_0$,
an effective integral Cartier divisor $D_0$ with support
$X_0\setminus M$, and an ample $\QQ$-line bundle $A$ on $X_0$.
For every pair of integers $a,j\geq1$,
\Cref{rem:number-field-approximation} gives a projective
modification $q_{a,j}:Y_{a,j}\to X_0$ over $K$, with $Y_{a,j}$
smooth and $q_{a,j}$ an isomorphism over $M$, and a semiample
$\QQ$-line bundle $N_{a,j}$, with a specified isomorphism
\[
 aG_{g,M}+A|_M\xrightarrow{\sim}N_{a,j}|_M,
\]
such that
\begin{equation}\label{eq:GS-approximation}
 \bnorm{N_{a,j}-(aG_g^{\geom}+A)}\leq j^{-1}.
\end{equation}
Here one takes $\varepsilon=1/a$ and $\delta=1/(aj)$
and multiplies the approximation by $a$.
The line bundle $A$ is fixed for all $a$ and $j$.

Nefness is preserved by algebraic pullback to smooth
quasi-projective varieties \cite[Section~2.5.5]{YZ}.
It is also preserved by extension of the ground field,
by base change of the approximations in \eqref{eq:GS-approximation}.
\end{remark}

\section*{Index of notation}
\phantomsection\label{sec:notation}
\addcontentsline{toc}{section}{Index of notation}
\begingroup
\small
\setlength{\tabcolsep}{4pt}
\renewcommand{\arraystretch}{1.12}
\begin{longtable}{@{}>{\raggedright\arraybackslash}p{.18\textwidth}>{\raggedright\arraybackslash}p{.50\textwidth}>{\raggedright\arraybackslash}p{.25\textwidth}@{}}
\toprule
Symbol & Meaning & Reference \\
\midrule
\endfirsthead
\toprule
Symbol & Meaning & Reference \\
\midrule
\endhead
\bottomrule
\endfoot
$U,X_0,D_0$ & Open variety, fixed compactification, and full Cartier boundary. & \Cref{sec:geometric-conventions} \\
$D_X=p^*D_0$ & Boundary divisor on a model $U\subset X$, pulled back along $p$. & \eqref{eq:full-boundary} \\
$L_U,\alpha_P$ & Underlying $\QQ$-line bundle and its specified isomorphism with a model restriction. & \Cref{def:model} \\
$\LL,P$ & Compactified line bundle and a model line bundle, also viewed as its constant class. & \Cref{def:completion} \\
$\AdPic(U/k)_{\QQ}$ & Group of compactified $\QQ$-line bundles over a field $k$. & \Cref{def:completion} \\
$\Delta(Q,P)$ & Boundary divisor determined by the specified comparison over $U$. & \eqref{eq:comparison-divisor} \\
$\|\cdot\|_{D_0}$ & Boundary norm on model differences and their completion. & \Cref{def:boundary-norm}; \eqref{eq:completed-order} \\
$h,\varphi_{h,e}$ & Metric and local weight $-\log\|e\|_h$. & \Cref{def:weights} \\
$\LL,h$ & Compactified line bundle with a compactified metric. & \Cref{def:compactified-metric} \\
$g_0,g_X$ & Boundary Green function and its pullback. & \Cref{def:boundary-green} \\
$t,P_t,h_t$ & Approximation tolerance, model line bundle, and model metric. & \eqref{eq:simultaneous-approximation} \\
$B_t,w_t$ & Boundary-corrected line bundle and its extended psh weight. & \Cref{prop:small-extension} \\
$A,A_t,\lambda$ & Fixed ample line bundle, its ample perturbation on a model, and the boundary tolerance. & \Cref{lem:ample-adjustment} \\
$\varepsilon,E_{t,\varepsilon}$ & Positive rational parameter and $B_t+\varepsilon A_t$. & \Cref{lem:exact-generation} \\
$\JJ(mw_t)$ & Coherent analytic multiplier ideal sheaf. & \Cref{def:multiplier-ideal} \\
$\II_{t,m}$ & Algebraic ideal corresponding to $\JJ(mw_t)$. & \Cref{prop:ideal-principalization} \\
$\rho_{t,m},F_{t,m}$ & Principalization and effective divisor of the image ideal. & \eqref{eq:principalized-ideal} \\
$\mathscr S_{t,\varepsilon,m}$ & Globally generated sheaf $\OO(mE_{t,\varepsilon})\otimes\II_{t,m}$. & \eqref{eq:section-sheaf} \\
$N_{t,\varepsilon,m}$ & Semiample model approximation on $Y_{t,m}$. & \eqref{eq:candidate-model} \\
$\overline{\LL},\LL^{\geom}$ & Arithmetic adelic line bundle and its geometric part. & \Cref{def:arithmetic-adelic,def:geometric-part} \\
$\sigma,h_\sigma$ & Complex embedding and the prescribed metric at that embedding. & \Cref{def:geometric-part} \\
$\overline G^{\,Y}_{\mathcal C/S}$ & Yuan's Gross--Schoen line bundle. & \Cref{def:yuan-GS} \\
$\overline{\OO}(\Phi_S)$ & Adelic line bundle associated with the relative $\varphi$-invariant. & \eqref{eq:phi-line} \\
$\overline G_g,G_g^{\geom}$ & Gao--Zhang's adelic Deligne pairing and its geometric part. & \Cref{def:GZ-GS} \\
$G_{g,M},\beta_{\GS}$ & Underlying Gross--Schoen line bundle and its complex curvature form. & \Cref{def:GZ-GS,prop:GS-metric} \\
$Q,N_{\CC},N_0$ & Reference line bundle, complex approximation, and approximation over the field of definition. & \Cref{prop:finite-transfer} \\
$\mathscr B_F,\mathscr B$ & Tensor product of the Galois conjugates and its descent; $N_0=\mathscr B/(br)$. & \eqref{eq:normalized-line} \\
\end{longtable}
\endgroup

\Needspace{8\baselineskip}
\begingroup
\raggedright
\phantomsection
\addcontentsline{toc}{section}{References}
\bibliographystyle{amsplain}
\bibliography{references}
\endgroup

\end{document}